\documentclass[12pt,leqno,a4paper]{article}
\usepackage{amsfonts}
\usepackage{amsmath, amssymb, amsthm, amsfonts}

\numberwithin{equation}{section}

\begin{document}
\title{On the  mean square of the error term for the two-dimensional divisor problems(II)}
\author{ Xiaodong Cao, Wenguang  Zhai  }
\date{}

\footnotetext[0]{2000 Mathematics Subject Classification: 11N37.}
\footnotetext[0]{Key Words: two-dimensional divisor problems , error
term, mean square, asymptotic formula .  } \footnotetext[0]{This
work is supported by National Natural Science Foundation of
China(Grant No. 10771127) .} \maketitle

{\bf Abstract.} Let $\Delta(a,b;x)$ denote the error term of the
general two-dimensional divisor problem. In this paper we shall
study the relation between the discrete mean value $\sum_{n\leq
T}\Delta^2(a,b;n)$ and the continuous mean value
$\int_1^T\Delta^2(a,b;x)dx$.

\section{\bf Introduction and state of results }

Suppose $1\leq a\leq b$ are two fixed integers. Without loss of
generality, we suppose  $(a,b)=1.$ Define
$$d(a,b;n):=\sum_{n=h^ar^b}1,\ \ \  \ D(a,b;x):=\sum_{n\leq x}
d(a,b;n).$$
 The two-dimensional divisor problems is to study the
error term
\begin{eqnarray}
\Delta(a,b;x):=  D(a,b;x)- \left\{\begin{array}{ll}
\left(\zeta(b/a)x^{1/a}+\zeta(a/b)x^{1/b}\right),&\mbox{if $1\leq a<b$,}\\
\left( x\log x+(2\gamma-1)x\right),& \mbox{if $a=b=1.$}
\end{array}\right.
\end{eqnarray}
  When $a=b=1,$ $\Delta(1,1;x)$ is the error term of the
well-known Dirichlet divisor problem. Dirichlet first proved that
 $\Delta(1,1;x)=O(x^{1/2}).$
 The exponent $1/2$ was improved by many authors.
The latest result reads

\begin{equation}
\Delta(1,1;x)\ll x^{131/416}\log^{26947/8320} x,
\end{equation}
which can be found in  Huxley\cite{Hu}. It is conjectured that
\begin{equation}
\Delta(1,1;x)=O(x^{1/4+\varepsilon}),
\end{equation}
which is supported by the classical mean-square  result
\begin{equation}
\int_1^T\Delta^2(1,1;x)dx=\frac{(\zeta(3/2))^4}{6\pi^2
\zeta(3)}T^{3/2} +O(T^{5/4+\varepsilon} )
\end{equation}
proved in \cite{Cr}. The estimate $O(T^{5/4+\varepsilon} )$ was
improved to $O(T\log^5 T)$ in \cite{T} and $O(T\log^4 T)$ in
\cite{P}. The mean square  of the error term  in (1.4) was studied
in  \cite{LT} and \cite{Ts1}.
 The higher-power moments of $\Delta(1,1;x)$ were
studied in\cite{H1,IS, Ts, Zw1, Zw2}.

When $a\not= b$, Richert \cite{Ri} proved that

\begin{eqnarray}
\Delta (a,b;x)\ll
\left\{
\begin{array}{ll}
x^{\frac {2}{3(a+b)}},&\mbox{if
$b\leq 2a$ }\\
x^{\frac {2}{5a+2b}},
&\mbox{if $b\geq 2a$},
\end{array}
\right.
\end{eqnarray}
Better upper estimates can be found in \cite{Kr1,Kr2,Ra,Sc}.
Hafner \cite{Ha2} showed that
 \begin{eqnarray}
 \Delta (a,b;x)= \Omega_{+}\left(x^{\frac {1}{2(a+b)}}(\log x)^{\frac {b}{2(a+b)}}\log \log
 x\right),
 \end{eqnarray}
and
\begin{eqnarray}
\Delta (a,b;x)= \Omega_{-}\left(x^{\frac {1}{2(a+b)}}e^{U(x)}\right),
\end{eqnarray}
where
\begin{eqnarray}
U(x)=B\left(\log \log x\right)^{\frac {b}{2(a+b)}}\left(\log \log \log x\right)^{\frac {b}{2(a+b)}-1},
\end{eqnarray}
for some constant $B>0$.

  For $1\leq a<b$ it is conjectured that the estimate
 \begin{eqnarray}
 \Delta (a,b;x)= O\left(x^{\frac {1}{2(a+b)}+\varepsilon}\right)
 \end{eqnarray}
 holds for $x\geq 2,$ which is supported partially by results of Ivi\'c\cite{I0}.
  Ivi\'c   showed that
\begin{eqnarray}
\int_1^T\Delta^2 (a,b;x)dx
\left\{
\begin{array}{ll}
\ll T^{1+\frac {1}{a+b}}\log^2T,\\
=\Omega(T^{1+\frac {1}{a+b}}).
\end{array}
\right.
\end{eqnarray}
The $\Omega$ result in (1.10) was improved by the first-named
author\cite{Cao} to
\begin{eqnarray}
  \int_1^T\Delta^2 (a,b;x)dx\gg T^{1+\frac {1}{a+b}}.
  \end{eqnarray}
Ivi\'c\cite{I0}
   conjectured that the
asymptotic formula
 \begin{eqnarray}
 \int_1^T\Delta^2 (a,b;x)dx =c_{a,b}(1+o(1)) T^{1+\frac {1}{a+b}}
 \end{eqnarray}
should hold for some constant $c_{a,b}>0$. This conjecture was
solved completely  in \cite{Zc}, where we proved   that if $1\leq a<
b$ and $(a,b)=1$ , then  for $T\geq 2$
\begin{eqnarray}
\int_1^T\Delta^2(a,b;x)dx=c_{a,b}T^{\frac{1+a+b}{a+b}}+O(T^{\frac{1+a+b}{a+b}-\frac{a}{2b(a+b)(a+b-1)}
}\log^{7/2} T),
\end{eqnarray}
where
$$c_{a,b}:=\frac{a^{b/(a+b)}b^{a/(a+b)}}{2(a+b+1)\pi^2}\sum_{n=1}^\infty
g_{a,b}^2(n)$$ and
$$g_{a,b}(n):=\sum_{n=h^ar^b}h^{-\frac{a+2b}{2a+2b}}r^{-\frac{b+2a}{2a+2b}}.$$

\bigskip

The aim of this paper is to study the relation between discrete mean
and continuous mean of $\Delta(a,b;x).$ This kind of problem is very
important and interesting in number theory.

  Vorono\"i\cite{Vo}
 essentially showed that for $x\geq 1$, the asymptotic formula
 \begin{eqnarray}
\sum_{n\leq x}\Delta(1,1;n)= \left(\frac{1}{2}-\psi (x)
\right)\Delta(1,1;x)+\int_1^x \Delta(1,1;t)dt
\end{eqnarray}
\begin{eqnarray*}
+\frac {1} {2}\left(\log x +2\gamma -1\right)x+O(\log x)
\end{eqnarray*}
holds.  For the mean square case, Hardy\cite{Har} proved that
$$\sum_{n\leq x}\Delta^2(1,1;n)=\int_1^x \Delta^2(1,1;t)dt+O(x^{1+\varepsilon}),$$
which was improved by  Furuya \cite{Fu} substantially to
\begin{eqnarray*}
\sum_{n\leq x}\Delta^2(1,1;n)&&=\int_1^x \Delta^2(1,1;t)dt+\frac 1 6 x\log^2x\\
&&+\frac {8\gamma -1}{12} x\log x+\frac {8\gamma^2-2\gamma+1}{12}x+
\left\{
\begin{array}{ll}
 O(x^{\frac 3 4}\log x),\\
\Omega_\pm (x^{\frac 3 4}\log x),
\end{array}
\right.
\end{eqnarray*}

For $1\leq a<b$, following  Vorono\"i\cite{Vo} we can get
\begin{eqnarray}
 \sum_{n\leq x}\Delta(a,b;n)= \left(\frac{1}{2}-\psi (x)
 \right)\Delta(a,b;x)+\int_1^x \Delta(a,b;t)dt
 \end{eqnarray}
 \begin{eqnarray*}
 +\frac {1} {2}\left(\zeta
  (\frac b a)x^{\frac 1 a}+\zeta (\frac a b)x^{\frac 1
  b}\right)+O(1).
 \end{eqnarray*}
From (1.14) and (1.15)    we may write
\begin{eqnarray*}
\sum_{n\leq x}\Delta(a,b;n)&&= \frac{1}{2}D(a,b;x)-\psi
(x)\Delta(a,b;x) + \int_1^x \Delta(a,b;t)dt\\&&\ \ \ \ +\left\{
\begin{array}{ll}
 O(\log x),&\mbox{if $a=b=1$},\\
O(1),&\mbox{if $1\leq a<b$},
\end{array}
\right.
\end{eqnarray*}
 which combining (1.5) and Lemma 2.2 implies that
 \begin{eqnarray}
 \sum_{n\leq x}\Delta(1,b;n)= \left(\frac 1 4+\frac {1} {2}\zeta
 (b)\right)x+
\left\{
\begin{array}{ll}
O\left(x^{1-\frac {1}{2(1+b)}}\right)\\
\Omega_{\pm}\left(x^{1-\frac {1}{2(1+b)}}\right)
\end{array}
\right.
\end{eqnarray}
for  $ b>1$,  and implies that
\begin{eqnarray}
 \sum_{n\leq x}\Delta(a,b;n)= \frac 1 4 x+
\left\{
\begin{array}{ll}
O\left(x^{1-\frac {1}{2(a+b)}}\right),\\
\Omega_{\pm}\left(x^{1-\frac {1}{2(a+b)}}\right),
\end{array}
\right.
\end{eqnarray}
for   $2\leq a<b.$ We omit the proofs of (1.15)-(1.17).

We shall study the mean square case for $1\leq a<b.$ Our main result
is the following

{\bf Theorem 1.} Let $1\leq a<b, (a,b)=1$,   and $x\geq 2$, then we
have
 \begin{eqnarray}
 \sum_{n\leq x}\Delta^2(a,b;n)=&& \left(\frac{1}{2}-\psi (x)
 \right)\Delta^2(a,b;x)+\int_1^x \Delta^2(a,b;t)dt
 \end{eqnarray}
\begin{eqnarray*}
+\frac {1} {4}\left(\zeta(\frac b a )x^{\frac 1 a}+\zeta(\frac a b)x^{\frac 1 b}\right)
+\left(\frac 1 a\zeta(\frac b a)x^{\frac 1 a-1}
+\frac 1 b\zeta(\frac a b)x^{\frac 1 b-1}\right) G_{(a,b)}(x)
\end{eqnarray*}

\begin{eqnarray*}
+\left\{
\begin{array}{ll}
\frac {1}{6}\zeta(b)\left(\zeta (b)x+2\zeta (\frac 1 b)x^{\frac 1 b}\right)+O(x^{1-\frac {3}{2(1+b)}}),&\mbox{if
$1=a<b,$}\\
O(x^{\frac 1 a-\frac {3}{2(a+b)}}), &\mbox{if $2\leq a<b,$}
\end{array}
\right.
\end{eqnarray*}
where
\begin{eqnarray}
\ \ \ \  G_{(a,b)}(x):=c_0x^{1-\frac {1}{2(a+b)}}
 \sum_{n=1}^{\infty}\frac {d^{*}(a,b;n)}{n^{1+\frac {1}{2(a+b)}}}
 \cos \left( 2(a+b)\pi (a^{-a}b^{-b}nx)^{\frac {1}{a+b}}+\theta _0\right) ,
 \end{eqnarray}
  \begin{eqnarray}
 &&d^{*}(a,b;n)=\sum_{m^an^b=n}m^{a-1}n^{b-1},\\
 &&c_0=\frac {(a^a
b^b)^{1+\frac{1}{2(a+b)}}}{2\pi^2\sqrt{ab(a+b)}},\ \
\theta_0=-\frac{ 3 \pi}{4}.
\end{eqnarray}

By Theorem 1 and (1.5), it is easy  to see that

{\bf Corollary 1.} Suppose $1\leq a<b, (a,b)=1$  and $x\geq 2$,
then we have
\begin{eqnarray}
 \sum_{n\leq x}\Delta^2(a,b;n)=\int_1^x \Delta^2(a,b;t)dt+
 R^*(a,b;x)
\end{eqnarray}
\begin{eqnarray*}
 +\left\{
\begin{array}{ll}
\frac {1}{12}\zeta(b)(3+2\zeta (b))x, &\mbox{if
$1=a<b,$}\\
\frac 1 4\zeta(\frac b a)x^{\frac 1 a}+\frac 1 4\zeta(\frac a
b)x^{\frac 1 b}, &\mbox{if $2\leq a<b,$}
\end{array}
\right.
\end{eqnarray*}
where $R^*(a,b;x)=O\left(x^{\frac 1 a -\frac {1}{2(a+b)}}\right)$
and $R^*(a,b;x)=\Omega_\pm\left(x^{\frac 1 a -\frac
{1}{2(a+b)}}\right)$.

{\bf Remark.} Generally speaking, the term $\frac 1 4\zeta(\frac a
b)x^{\frac 1 b}$ in Corollary 1 can not be removed since $\frac 1
b>\frac 1 a -\frac {1}{2(a+b)}$ for $b<\frac {1+\sqrt{17}}{4}a$.

From Corollary 1 and (1.13) we get

 {\bf Corollary 2.} Suppose $1\leq a<b, (a,b)=1$ and $x\geq 2$, then
 \begin{eqnarray*}
 \sum_{n\leq x}\Delta^2(a,b;n)=
 c_{a,b}x^{\frac{1+a+b}{a+b}}+O(x^{\frac{1+a+b}{a+b}-\frac{a}{2b(a+b)(a+b-1)}
 }\log^{7/2} x).
 \end{eqnarray*}

We have also the following Theorem 2, which slightly improves
Furuya's result.

{\bf Theorem 2.} For $x\geq 2$, we have
 \begin{eqnarray}
 \sum_{n\leq x}\Delta^2(1,1;n)=&& \left(\frac{1}{2}-\psi (x)
 \right)\Delta^2(1,1;x)+\int_1^x \Delta^2(1,1;t)dt
 \end{eqnarray}
\begin{eqnarray*}
+\frac 1 6 x\log^2x+\frac {8\gamma -1}{12} x\log x+\frac {8\gamma^2-2\gamma+1}{12}x
 +(\log x+2\gamma)G_{(1,1)}(x)+O(\sqrt{x}\log x),
\end{eqnarray*}
where
\begin{eqnarray*}
 G_{(1,1)}(x):=\frac {1}{2{\sqrt 2} \pi^2}x^{\frac 3 4}
 \sum_{n=1}^{\infty}\frac {d(n)}{n^{\frac 5 4}}
 \sin \left( 4\pi \sqrt {nx}-\frac {\pi}{4}\right).
 \end{eqnarray*}

Our proof is based on the method of  Furuya \cite{Fu}. We need   a
sharper asymptotic formula for the error term $\Delta(a,b;x)$, and
then
  evaluate a kind of integrals involving the $\psi$-function.

{\bf Notations.} For a real number $u,$ $[u]$ denotes the greatest
integer not exceeding $u$, $\psi(u)=u-[u]-1/2$. Let $(m,n)$ denote
the greatest common divisor of natural numbers $m$ and $n$. $n\sim
N$ means $N<n\leq 2N.$ $\varepsilon$ always denotes a sufficiently
small positive constant. In this paper, the constants implied by
$O$ depend olny on $a, b$ or  $\varepsilon$ when it occures.
$f(x)=O(g(x))$ or $f(x)\ll g(x)$ means that $|f(x)|\leq Cg(x)$ for
$x\geq x_0$ and some absolute constant $C>0$. $f(x)=\Omega_\pm
(g(x))$ means that both $f(x)=\Omega_+ (g(x))$ and $f(x)=\Omega_-
(g(x))$ holds.

\section{\bf Some preliminary Lemmas}

In order to prove our theorems, we need the following  lemmas.

{\bf Lemma 2.1.}  Let $f(n)$ be an arithmetic function, and $E(x)$
be  the   error term defined by
\begin{eqnarray}
E(x):=\sum_{n\leq x}f(n) -g(x).
\end{eqnarray}
Suppose  $g(x)$ is continuously differentiable. For any fixed
positive integer $k$, we have
\begin{eqnarray}
\sum_{n\leq x}E^{k} (n)= &&\left(\frac{1}{2}-\psi (x)
\right)E^k(x)+\int_1^x E^k (u)du
\end{eqnarray}
\begin{eqnarray*}
&&+k \int_1^x \left(\frac{1}{2}-\psi (u) \right )
g'(u)E^{k-1}(u)du.
\end{eqnarray*}

{\bf Proof.} This is Lemma 1 of Furuya \cite{Fu}.

{\bf Lemma 2.2.} (Vorono\"i type formula) Let $x\geq 1$, and
$d^{*}(a,b;n)$ is defined by (1.20), then for every fixed positive
integer $q$,  the following asymptotic formula holds
\begin{eqnarray}
\int_0^x \Delta (a,b;t)dt=\frac x 4+\zeta (-a)
 \zeta (-b)+O(x^{1-\frac {1}{2(a+b)}-\frac {q}{a+b}})+\end{eqnarray}
\begin{eqnarray*}
 \sum_{m=0}^{q-1}c_m \left(
 \sum_{n=1}^{\infty}\frac {d^{*}(a,b;n)}{n^{1+\frac {1}{2(a+b)}+\frac {m}{a+b}}}
 \cos \left( 2(a+b)\pi (a^{-a}b^{-b}nx)^{\frac {1}{a+b}}+\theta _m\right)\right) x^{1-\frac {1}{2(a+b)}-\frac {m}{a+b}}
 ,
 \end{eqnarray*}
where $c_m$ and $\theta_m$ are real numbers, $c_0$ and $ \theta_0$
are defined in Theorem 1.

{\bf Proof.} This is Theorem 3 of the first-named author \cite{Cao}.
We note that when $a=b=1$, Lemma 2.2 was already proved by Tong
\cite{T}.

Now for $x\geq 1$, we define the $\psi_j (x)$ by the following
recurrence relation
\begin{eqnarray}
\psi_j (x):=\int_1^x \psi_{j-1}(t)dt,  (j=1,2,\cdots),
\end{eqnarray}
 for convenience, we  use the notation $\psi_0 (x)=\psi
(x)$ here.

{\bf Lemma 2.3.} For $x\geq 1$, we have
\begin{eqnarray}
\int_1^x \psi ^{2k-1}(t)dt=\frac {1}{2k}\left(\psi ^{2k}(x)-\frac
{1}{2^{2k}}\right), (k=1,2,\cdots)
\end{eqnarray}
and
\begin{eqnarray}
\int_1^x \psi ^{2k}(t)dt=\frac {1}{2k+1}\left(\frac
{x-1}{2^{2k}}-\frac {1}{2^{2k}}\psi(x)+\psi ^{2k+1}(x)\right),
(k=1,2,\cdots).
\end{eqnarray}
{\bf Proof.} We prove (2.6) only, the proof of (2.5) is similar. By
  simple calculations, we get
\begin{eqnarray*}
\int_1^x \psi ^{2k}(t)dt&&=\sum_{n=0}^{[x]-1}\int_{n}^{n+1} \psi
^{2k}(t)dt+\int_{[x]}^x \psi ^{2k}(t)dt\\
&&=([x]-1)\int_1^2 \psi ^{2k}(t)dt+\int_{[x]}^x (t-[x]-\frac 1 2) ^{2k}dt\\
&&=2([x]-1)\int_0^{\frac {1}{2}} u ^{2k}du+\int_{-\frac {1}{2}}^{\psi (x)} u^{2k}du\\
&&=([x]-1)\frac {1}{2^{2k}(2k+1)}+\frac {1}{(2k+1)}\left( \psi
^{2k+1}(x)+\frac {1}{2^{2k+1}}\right),
 \end{eqnarray*}
and whence (2.6) follows.

From Lemma 2.3, we  easily get that for $x\geq 1$
\begin{eqnarray}
\psi_1(x)=\frac 1 2(\psi^2(x)-\frac 1 4),
\end{eqnarray}
and
\begin{eqnarray}
\psi_2(x)=-\frac {1} {12}x+\frac 16 \psi^3(x)-\frac
{1}{24}\psi(x)+\frac {1}{12}.
\end{eqnarray}

{\bf Lemma 2.4.} Let $x\geq 1$ and define $W_\alpha(x):=\int_1^x
t^\alpha \psi (t)dt$. If $\alpha \not = -1, -2$, we have

\begin{eqnarray}
W_\alpha (x)=-\frac {x^{\alpha +2}}{(\alpha +1)(\alpha +2)}+\frac
{1}{\alpha +1}\sum_{n\leq x} n^{\alpha +1} +\end{eqnarray}
\begin{eqnarray*}
\frac {\psi (x)}{\alpha +1}x^{\alpha +1}-\frac {\alpha}{2(\alpha
+1)(\alpha +2)},
\end{eqnarray*}

\begin{eqnarray}
W_{-2} (x)&&=\log x -\psi (x)x^{-1}-\sum_{n\leq x}n^{-1}+\frac{1}{2}\\
&&=\frac{1}{2}-\gamma +\left(\psi_1(x)+\frac
{1}{12}\right)x^{-2}+O(x^{-3}),
\end{eqnarray}
and
\begin{eqnarray}
W_{-1} (x)=2\int_{1}^{\infty}\frac
{\psi_2(t)}{t^3}dt+\left(\psi_1(x)+\frac
{1}{12}\right)x^{-1}+O(x^{-2}).
\end{eqnarray}

{\bf Proof.}  First, we suppose $x\geq 2$. Similar to  the proof of
Lemma 2.3, we have
\begin{eqnarray}
W_\alpha (x)=\sum_{n=1}^{[x]-1}\int_{n}^{n+1} t^\alpha (t-n-\frac
 1 2)dt +\int_{[x]}^x t^\alpha (t-[x]-\frac 1 2)dt
\end{eqnarray}
\begin{eqnarray*}
=&&\left(\sum_{n=1}^{[x]-1}\int_{n}^{n+1} t^{\alpha +1}dt
+\int_{[x]}^x t^{\alpha +1}dt \right)
-\left(\sum_{n=1}^{[x]-1}(n+\frac 1 2)\int_{n}^{n+1} t^{\alpha}dt +([x]+\frac 1 2)\int _{[x]}^x t^{\alpha }dt\right)\\
=&&\int_{1}^{x} t^{\alpha +1}dt-\frac {1}{\alpha
+1}\left(\sum_{n=1}^{[x]-1}(n+\frac 1 2)\left((n+1)^{\alpha +1}-n^{\alpha +1}\right)
+([x]+\frac 1 2)(x^{\alpha +1}-[x]^{\alpha +1})\right)\\
=&&\frac {x^{\alpha +2}-1}{\alpha +2}-\frac {1}{\alpha
+1}\left(\sum_{n=2}^{[x]}(n-\frac 1 2)n^{\alpha
+1}-\sum_{n=1}^{[x]-1}(n+\frac 1 2)n^{\alpha +1}\right) -\frac
{1}{\alpha +1}([x]+\frac 1 2)(x^{\alpha +1}-[x]^{\alpha +1})\\
=&&\frac {x^{\alpha +2}-1}{\alpha +2}-\frac {1}{\alpha
+1}\left(-\sum_{n=2}^{[x]-1}n^{\alpha +1}+([x]-\frac 1
2)[x]^{\alpha +1}-\frac 3 2\right)-\frac {1}{\alpha +1}([x]+\frac
1 2)(x^{\alpha +1}-[x]^{\alpha +1}).
\end{eqnarray*}
We make some simplification, and obtain that (2.9) holds in this
case. Next, if $1\leq x<2$, then $[x]=1$, it is easy to check that
(2.9) also holds, and this completes the proof of (2.9).

Now we consider the case $\alpha =-2$. For  $x\geq 2$, by the same method as
above,  we have
\begin{eqnarray}
W_{-2} (x)=\log x -\psi (x)x^{-1}-\sum_{n\leq
x}n^{-1}+\frac{1}{2}.
\end{eqnarray}
If $1\leq x<2$, then $W_{-2} (x)=\log x -\frac 3 2+\frac 3 2
x^{-1}$, and (2.14) also holds in this case. This completes the
proof of (2.10). Applying Euler-Maclaurin formula(see (2.20)) to
the sum $\sum_{n\leq x}n^{-1}$, we can get (2.11).

Finally, by applying integration by parts and (2.8), we have
\begin{eqnarray*}
W_{-1} (x)=&&\psi_1 (x)x^{-1}+\int_1^x t^{-2} \psi_1 (t)dt\\
=&&\psi_1 (x)x^{-1}+\psi_2 (x)x^{-2}+2\int_1^x t^{-3} \psi_2
(t)dt\\
=&&2\int_1^{\infty} t^{-3} \psi_2 (t)dt+\psi_1 (x)x^{-1}+\psi_2
(x)x^{-2}-2\int_x^{\infty} t^{-3} \psi_2 (t)dt\\
=&&\int_1^{\infty} t^{-3} \psi_2 (t)dt+\psi_1 (x)x^{-1}+(-\frac
{x}{12}+O(1))x^{-2}-2\int_x^{\infty} \frac {-\frac
{t}{12}+O(1)}{t^{3} }(t)dt.
\end{eqnarray*}
This finishes the proof of Lemma 2.4.

 When $\alpha =0$,  from (2.9) and some easy calculations, we get for  $x\geq 1$
that
\begin{eqnarray}
W_0(x)=\psi_1(x)=\frac{1}{2} ( \psi^2(x)-\frac{ 1}{4}),
\end{eqnarray}
and
\begin{eqnarray}
W_0(n)=\psi_1(n)=0.\end{eqnarray}

   If $\alpha =1$,  we get for $x\geq 1$ that
\begin{eqnarray}
W_1(x)=-\frac {x^3}{6}+\frac 1 2\frac {[x]([x]+1)(2[x]+1)}{6}+
\frac {\psi (x)}{2}x^2-\frac {1}{12}
\end{eqnarray}
\begin{eqnarray*}
 =-\frac {1}{24}x+\frac 1 2\psi^2 (x)x-\frac 1 6\psi^3(x)+\frac
 {1}{24}\psi(x)-\frac {1}{12}.
\end{eqnarray*}

When $\alpha$ is a non-negative integer, we may use the well-known
Bernoulli polynomial to express $W_\alpha(x)$. Otherwise, we can use
the following Lemma 2.5 to estimate it.

{\bf Lemma 2.5.} Let $\alpha \not = -1, -2$, then for $x\geq 1$,
 \begin{eqnarray}
 W_\alpha(x)=\frac {1}{\alpha +1}\left(\zeta (-1-\alpha)-\frac {\alpha}{2(2+\alpha)}\right)+(\psi_1(x)
 +\frac {1}{12})x^\alpha +O(x^{\alpha -1})
\end{eqnarray}

{\bf Proof.} From Euler-Maclaurin formula, for $s\not =1$ and
$x\geq 1$, we have
\begin{eqnarray}
\sum_{n\leq x}n^{-s}= \zeta (s)+\frac{x^{1-s}}{1-s}-\psi
(x)x^{-s}-s(\psi_1(x)+\frac {1}{12})x^{-s-1}+O(x^{-s-2}),
\end{eqnarray}
and
\begin{eqnarray}
\sum_{n\leq x}n^{-1}=\log x+\gamma-\psi (x)x^{-1}-(\psi_1(x)+\frac
{1}{12})x^{-2}+O(x^{-3}).
\end{eqnarray}
Now (2.18) is a immediate consequence of (2.9) and (2.19).

\section{\bf An asymptotic formula for the error term  $\Delta(a,b; x)$}
It is well-known that
\begin{eqnarray}
\Delta (a,b;x)=-\sum_{n^{a+b}\leq x}
\psi\left(\left(\frac{x}{n^{b}}\right)^{\frac
{1}{a}}\right)+\psi\left(\left(\frac{x}{n^{a}}\right)^{\frac {1}{
b}}\right)+O(1).
\end{eqnarray}
However, the error function $O(1)$ in (3.1) is too large to prove
our theorems. So we need a sharper form than (3.1). In this section
 we shall prove such a   lemma.

{\bf Lemma 3.1.}  Let $(a,b)=1$ and $x\geq 1$, we define the error
function $R(a,b;x)$ by
\begin{eqnarray}
R(a,b;x):=\Delta (a,b;x)+\sum_{n^{a+b}\leq x}
\psi\left(\left(\frac{x}{n^{b}}\right)^{\frac
{1}{a}}\right)+\psi\left(\left(\frac{x}{n^{a}}\right)^{\frac {1}{
b}}\right).
\end{eqnarray}
 Then one has
\begin{eqnarray}
R(a,b;x)=-\frac {(a+b)^2}{ab}\psi_1(x^{\frac {1}{a+b}})+\frac
{b(a+b)}{a^2}x^{\frac 1 a}\int_{x^{\frac
{1}{a+b}}}^{\infty}\psi_1(t)t^{-2-\frac b a}dt
\end{eqnarray}
\begin{eqnarray*}
+\frac {a(a+b)}{b^2}x^{\frac 1 b}\int_{x^{\frac
{1}{a+b}}}^{\infty}\psi_1(t)t^{-2-\frac a b}dt.
\end{eqnarray*}
In particular, we have
\begin{eqnarray}
R(a,b;x)=-\frac {(a+b)^2}{ab}\psi _1(x^{\frac {1}{a+b}})-\frac
{1}{12}\frac {a^2+b^2}{ab}+O\left(x^{-\frac {1}{a+b}}\right).
\end{eqnarray}
Furthermore, if $x^{\frac {1}{a+b}}$ is not an integer, then the
derivative of $R(a,b;x)$ satisfies
\begin{eqnarray}
R'(a,b;x)=&&-\frac {a+b}{ab}\psi(x^{\frac {1}{a+b}})x^{\frac
{1}{a+b}-1}-\frac {a^2+b^2}{ab}\psi_1(x^{\frac {1}{a+b}})x^{-1}
\end{eqnarray}
\begin{eqnarray*}
-\frac {a+b}{12} x^{- 1 }+O\left(x^{-1-\frac {1}{a+b}}\right).
\end{eqnarray*}

{\bf Proof.}  By applying the Dirichlet hyperbola method, we
easily obtain
\begin{eqnarray}
D(a,b;x)=\sum_{n\leq x}d(a,b;n)=\sum_{m^a n^b\leq x}1
\end{eqnarray}
\begin{eqnarray*}
=&&\sum_{m^{a+b}\leq x}\left[\left(\frac {x}{m^a}\right)^{\frac 1
b}\right]+\sum_{n^{a+b}\leq x}\left[\left(\frac
{x}{n^b}\right)^{\frac 1 a}\right]-\left(\sum_{m^{a+b}\leq
x}1\right)\left(\sum_{n^{a+b}\leq x}1\right)\\
=&&-\sum_{m^{a+b}\leq x}
\psi\left(\left(\frac {x}{m^a}\right)^{\frac 1 b}\right)+\psi\left(\left(\frac {x}{n^b}\right)^{\frac 1 a}\right)\\
&&+x^{\frac 1 b}\sum_{m^{a+b}\leq x}\frac {1}{m^{\frac a b
}}+x^{\frac 1 a}\sum_{n^{a+b}\leq x}\frac {1}{n^{\frac b a
}}-\left[x^{\frac {1}{a+b}}\right]^2-\left[x^{\frac
{1}{a+b}}\right].
\end{eqnarray*}

It is easy to see that the function $\psi_1(x)$ is a periodic
function  with period 1 and is therefore continuous. For real $x\geq
1$ and $s>0$, by using Riemann -Stieltjes integration, and then
integration by parts, we get that
\begin{eqnarray*}
\sum_{n\leq x}n^{-s}=&&\int_{1-0}^xt^{-s}d[t]=[x]x^{-s}+s\int_1^x [t] t^{-s-1} dt\\
=&&-\psi (x)x^{-s}+(x-\frac 1 2)x^{-s}+s\int_1^x (t-\frac 1 2)
t^{-s-1} dt\\
 &&-s\left(\psi_1 (x)x^{-s-1}+(s+1)\int_1^x t^{-s-2}
 \psi_1(t)dt\right).
\end{eqnarray*}
Hence, for $s>0$
\begin{eqnarray}
\sum_{n\leq x}n^{-s}=\left\{ \begin{array}{ll}
\frac{2x^{1-s}-1-s}{2(1-s)}-\frac {\psi (x)}{x^{s}}-\frac
{s\psi_1(x)}{x^{s+1}}-s(s+1)\int_1^x \frac {\psi_1(t)}{t^{s+2}}
 dt,&\mbox{if  $s\not =1$}\\
\log x+\frac 1 2-\frac {\psi (x)}{x}-\frac
{\psi_1(x)}{x^{2}}-2\int_1^x \frac {\psi_1(t)}{t^{3}}dt,&\mbox{if
$s=1$}.
\end{array}
\right.
\end{eqnarray}

It is well-known that for $s>0$ and $s\not =1$
\begin{eqnarray}
\zeta(s)=\frac {s+1}{2(s-1)}-s\int_1^\infty t^{-s-1} \psi (t)dt.
\end{eqnarray}
In addition, from (2.11) we have
\begin{eqnarray}\int_1^\infty t^{-2} \psi (t)dt=\frac 1 2 -\gamma.
\end{eqnarray}
Integrating by parts again, we see that (3.8) and (3.9) are
equivalent to
\begin{eqnarray}
\int_1^\infty t^{-s-2} \psi_1 (t)dt=\frac
{1}{2s(s-1)}-\frac{\zeta(s)}{s(s+1)},\mbox{ if $s>0$ and $s\not
=1$},
\end{eqnarray}
and
\begin{eqnarray}
\int_1^\infty t^{-3} \psi_1 (t)dt=\frac {1}{2}\left(\frac 1
2-\gamma\right),
\end{eqnarray}
respectively.
 Inserting (3.10) and (3.11) into (3.7), we obtain for $s>0$
\begin{eqnarray}
\ \ \ \ \ \sum_{n\leq x}n^{-s}=\left\{ \begin{array}{ll}
\frac{x^{1-s}}{1-s}+\zeta (s)-\frac {\psi (x)}{x^{s}}-\frac
{s\psi_1(x)}{x^{s+1}}+s(s+1)\int_x^\infty \frac {\psi_1(t)}{t^{s+2}}
 dt,&\mbox{if  $s\not =1$}\\
\log x+\gamma-\frac {\psi (x)}{x}-\frac
{\psi_1(x)}{x^{2}}+2\int_x^\infty \frac
{\psi_1(t)}{t^{3}}dt,&\mbox{if $s=1$}.
\end{array}
\right.
\end{eqnarray}

Now, taking $s=\frac b a$ and $s=\frac a b$ in (3.12) respectively,
then combining (1.1), (2.7 ), (3.2),  the following simple relations
\begin{eqnarray*}
\left[ x^{\frac{1}{a+b}}\right] =x^{\frac {1}{a+b}}-\psi\left(
x^{\frac {1}{a+b}}\right)-\frac {1}{2}
\end{eqnarray*}
and
\begin{eqnarray*}
 \left[x^{\frac{1}{a+b}}\right] ^2= x^{\frac {2}{a+b}}+\psi^2\left(
x^{\frac {1}{a+b}}\right)+\frac {1}{4}-2\psi\left( x^{\frac
{1}{a+b}}\right) x^{\frac {1}{a+b}}-x^{\frac {1}{a+b}}+\psi\left(
x^{\frac{1}{a+b}}\right),
\end{eqnarray*}
we can get (3.3).

   From formula (2.8) , it is easy to check that for $y\geq 1$
\begin{eqnarray*}
\psi_2(y)=\int_1^y\psi_1(t)dt=-\frac {y}{12}+O(1).
\end{eqnarray*}

  Integrating by parts, we obtain
\begin{eqnarray*}
&&\int_{x^{\frac {1}{a+b}}}^\infty \psi_1(t)t^{-\frac {b}{a} -2} dt\\
&&=-\left(x^{\frac{1}{a+b}}\right)^{-\frac {b}{a}
-2}\psi_2\left(x^{\frac{1}{a+b}}\right)+(2+\frac {b}{
a})\int_{x^{\frac {1}{a+b}}}^\infty \psi_2(t)t^{-\frac {b}{a} -3} dt\\
&&=-x^{-\frac {2a+b}{a(a+b)}}\left(-\frac {1}{12}x^{\frac
{1}{a+b}}+O(1)\right)+(2+\frac {b}{a})\int_{x^{\frac
{1}{a+b}}}^\infty
\left(-\frac {t}{12}+O(1)\right)t^{-\frac{b}{a} -3} dt\\
&&=- \frac{a}{12(a+b)} x^{-\frac {1}{a}}+O\left( x^{-\frac {1}{
a}-\frac {1}{a+b}}\right)  .
\end{eqnarray*}
Similarly, we also have
\begin{eqnarray*}
\int_{x^{\frac {1}{a+b}}}^\infty \psi_1(t)t^{-\frac {a}{b} -2}
dt=- \frac{b}{12(a+b)} x^{-\frac {1}{b}}+O\left( x^{-\frac {1}{
b}-\frac {1}{a+b}}\right)  .
\end{eqnarray*}
Combining the above two estimates and (3.3) completes the proof of
(3.4).

Finally, we suppose that $x^{\frac {1}{a+b}}$ is not an integer,
thus $\psi_1\left( x^{\frac {1}{a+b}}\right)$ is differentiable. By
differentiating the both sides of (3.3) with respect to $x$, and
then applying the above two estimates, we have
\begin{eqnarray*}
 R'(a,b;x)=&&-\frac {(a+b)^2}{ab}\frac {1}{a+b}\psi(x^{\frac
{1}{a+b}})x^{\frac {1}{a+b}-1}+\frac {b(a+b)}{a}x^{\frac 1 a
-1}\int_{x^{\frac {1}{a+b}}}^\infty \psi_1(t)t^{-\frac {b}{a} -2}
dt\\
&&-\frac {b(a+b)}{a}x^{\frac 1 a }\psi_1(x^{\frac
{1}{a+b}})x^{\frac {1}{a+b}(-\frac b a-2)}\frac {1}{a+b}x^{\frac
{1}{a+b}-1}+
\frac {a(a+b)}{b}x^{\frac 1 b -1}\int_{x^{\frac
{1}{a+b}}}^\infty \psi_1(t)t^{-\frac {a}{b} -2} dt\\
&& -\frac {a(a+b)}{b}x^{\frac 1 b}\psi_1(x^{\frac
{1}{a+b}})x^{\frac {1}{a+b}(-\frac a b-2)}\frac {1}{a+b}x^{\frac
{1}{a+b}-1}
\\
=&&-\frac {a+b}{a b}\psi (x^{\frac {1}{a+b}})x^{\frac
{1}{a+b}-1}-\frac {a^2+b^2}{a b}\psi_1(x^{\frac {1}{a+b}})x^{-1}-
\frac{a+b}{12} x^{- 1 }+O\left(x^{-1-\frac {1}{a+b}}\right) ,
\end{eqnarray*}
and this completes the proof of Lemma 3.1.

\section{\bf Integral formulas involving the $\psi$-function}

In this section we shall evaluate integrals involving the $\psi$
functions, which are important for our proof.

{\bf Lemma 4.1.} Let $n$ be an positive integer, $\alpha$ real,
$x\geq 1$, and $n^{a+b}\leq x$, we define
$I_{(a,b)}^{(\alpha)}(n,x):=\int _{n^{a+b}}^{x}t^\alpha \psi
(t)\psi\left(\left(\frac {t}{n^b}\right)^{\frac {1}{a}}\right)dt$.
Then
\begin{eqnarray}
I_{(a,b)}^{(\alpha)}(n,x)=\frac  {1}{n^{\frac {b}{a}}}\left(
W_{\alpha +\frac 1 a}(x)-W_{\alpha +\frac 1
a}(n^{a+b})\right)+\sum_{j=n}^{\left[\left(\frac
{x}{n^b}\right)^{\frac {1}{a}}\right]}W_{\alpha}(j^a
n^b)+\end{eqnarray}
 \begin{eqnarray*}
 (n-\frac 1 2)W_{\alpha}(n^{a+b})-\left( \left(\frac {x}{n^b}\right)^{\frac
{1}{a}}-\psi\left( \left(\frac {x}{n^b}\right)^{\frac {1}{a}}
\right)\right) W_{\alpha}(x),
\end{eqnarray*}
where $W_{\alpha}(x)$ is defined by Lemma 2.4.

 {\bf Proof.} We first suppose $\left[\left(\frac
{x}{n^b}\right)^{\frac {1}{a}}\right]\geq n+1$. We divide the
integral interval into the following subintervals, and obtain

\begin{eqnarray}
I_{(a,b)}^{(\alpha)}(n,x)=\sum_{j=n}^{\left[\left(\frac
{x}{n^b}\right)^{\frac {1}{a}}\right]-1}\int_{j^an^b}^{(j+1)^an^b}
t^\alpha \psi (t)\psi\left(\left(\frac {t}{n^b}\right)^{\frac
{1}{a}}\right)dt
\end{eqnarray}
\begin{eqnarray*}
&&+\int_{\left[\left(\frac {x}{n^b}\right)^{\frac
{1}{a}}\right]^an^b}^{x} t^\alpha \psi (t)\psi\left(\left(\frac
{t}{n^b}\right)^{\frac {1}{a}}\right)dt\\
=&&\sum_{j=n}^{\left[\left(\frac {x}{n^b}\right)^{\frac
{1}{a}}\right]-1}\int_{j^an^b}^{(j+1)^an^b} t^\alpha \psi
(t)\left(\left(\frac {t}{n^b}\right)^{\frac {1}{a}}-j-\frac 1
2\right)dt\\
&&+\int_{\left[\left(\frac {x}{n^b}\right)^{\frac
{1}{a}}\right]^an^b}^{x} t^\alpha \psi (t)\left(\left(\frac
{t}{n^b}\right)^{\frac {1}{a}}-\left[\left(\frac
{x}{n^b}\right)^{\frac {1}{a}}\right]-\frac 1 2\right)dt\\
=&&\frac {1}{n^{\frac b a}}\left(\sum_{j=n}^{\left[\left(\frac
{x}{n^b}\right)^{\frac {1}{a}}\right]-1}\int_{j^an^b}^{(j+1)^an^b}
t^{\alpha +\frac 1 a} \psi (t)dt+ \int_{\left[\left(\frac
{x}{n^b}\right)^{\frac {1}{a}}\right]^an^b}^{x}t^{\alpha +\frac 1
a} \psi (t)dt\right)\\
&&-\sum_{j=n}^{\left[\left(\frac {x}{n^b}\right)^{\frac
{1}{a}}\right]-1}(j+\frac 1 2)\int_{j^an^b}^{(j+1)^an^b} t^{\alpha
} \psi (t)dt -\left(\left[\left(\frac {x}{n^b}\right)^{\frac
{1}{a}}\right]+\frac 1 2\right)\int_{\left[\left(\frac
{x}{n^b}\right)^{\frac {1}{a}}\right]^an^b}^{x}t^{\alpha } \psi
(t)dt\\
=&&\frac {1}{n^{\frac b a}}\int_{n^{a+b}}^{x} t^{\alpha +\frac 1
a} \psi (t)dt-\sum_{j=n}^{\left[\left(\frac {x}{n^b}\right)^{\frac
{1}{a}}\right]-1}(j+\frac 1 2)\left(W_{\alpha}((j+1)^a
n^b)-W_{\alpha}(j^an^b)\right)\\
&& -\left(\left[\left(\frac {x}{n^b}\right)^{\frac
{1}{a}}\right]+\frac 1
2\right)\left(W_{\alpha}(x)-W_{\alpha}\left( \left[\left(\frac
{x}{n^b}\right)^{\frac {1}{a}}\right]^an^b\right)\right).
\end{eqnarray*}
Moreover, by Abel's summation formula, we get
\begin{eqnarray}
\sum_{j=n}^{\left[\left(\frac {x}{n^b}\right)^{\frac
{1}{a}}\right]-1}(j+\frac 1 2)\left(W_{\alpha}((j+1)^a
n^b)-W_{\alpha}(j^an^b)\right)
\end{eqnarray}
\begin{eqnarray*}
=&&\sum_{j=n+1}^{\left[\left(\frac {x}{n^b}\right)^{\frac
{1}{a}}\right]}(j-\frac 1 2)W_{\alpha}(j^a
n^b)-\sum_{j=n}^{\left[\left(\frac {x}{n^b}\right)^{\frac
{1}{a}}\right]-1}(j+\frac 1 2)W_{\alpha}(j^a n^b)\\
=&&-\sum_{j=n+1}^{\left[\left(\frac {x}{n^b}\right)^{\frac
{1}{a}}\right]-1}W_{\alpha}(j^a n^b)+\left(\left[\left(\frac
{x}{n^b}\right)^{\frac {1}{a}}\right]-\frac 1
2\right)W_{\alpha}\left( \left[\left(\frac {x}{n^b}\right)^{\frac
{1}{a}}\right]^an^b\right)-(n+\frac 1 2)W_\alpha (n^{a+b}).
\end{eqnarray*}
Combining (4.2) and (4.3), we find that (4.1) holds in this case.
If $\left[\left(\frac {x}{n^b}\right)^{\frac {1}{a}}\right]= n$,
it is easy to check that (4.1) also holds, and this completes the
proof of Lemma 4.1.

 {\bf Lemma 4.2.} Let $\alpha \not = -1, -2$, $n^{a+b}\leq
x$, and  $x\geq 1$. Then
\begin{eqnarray}
I_{(a,b)}^{(0)}(n,x)=\frac  {1}{n^{\frac {b}{a}}}\left( W_{\frac 1
a}(x)-W_{\frac 1 a}(n^{a+b})\right)
 -\left( \left(\frac {x}{n^b}\right)^{\frac
{1}{a}}-\psi\left( \left(\frac {x}{n^b}\right)^{\frac {1}{a}}
\right)\right) \psi_1(x).
\end{eqnarray}
If $\alpha\not=-\frac 1 a$, then
\begin{eqnarray}
I_{(a,b)}^{(\alpha)}(n,x)=\frac {1} {12(1+a\alpha)}\left(
\frac{x^{\alpha +\frac 1 a}}{n^{\frac b
a}}-n^{(a+b)\alpha+1}\right) + \psi\left( \left(\frac
{x}{n^b}\right)^{\frac {1}{a}} \right)\psi_1(x)x^\alpha
\end{eqnarray}
\begin{eqnarray*}
+O\left( \frac{x^{\alpha +\frac 1 a-1}}{n^{\frac b a}}\log
x+x^{\alpha -\frac 1 a}n^{\frac b a}+n^{(a+b)\alpha-1}
+n^{(a+b)(\alpha-1)+1}\log x\right),
\end{eqnarray*}

and if $\alpha =-\frac 1 a$, then
\begin{eqnarray}
I_{(a,b)}^{(-\frac 1 a)}(n,x)=\frac {1} {12 a}\frac {1}{n^{\frac b
a}}\left( \log x -(a+b)\log n\right) + \psi\left( \left(\frac
{x}{n^b}\right)^{\frac {1}{a}} \right)\psi_1(x)x^{-\frac 1 a}
\end{eqnarray}
\begin{eqnarray*}
+O\left( \frac{x^{-1}}{n^{\frac b a}}\log x+x^{-\frac 2 a}n^{\frac
b a}+n^{-\frac b a-2}+n^{-(a+b)(\frac 1 a+1)+1}\log x \right).
\end{eqnarray*}

{\bf Remark 2.} In fact, the Lemma 3 of Furuya \cite{Fu} is a
special case of (4.4) with $a=1$. In the present paper, when
$\alpha=0$, it is sufficient
 for us to apply a weaker estimate (4.5).

{\bf Proof.} First, (4.4) is an immediate consequence of Lemma 4.1
and (2.16).

By Lemma 4.1, Lemma 2.5 and some simplification, we have

\begin{eqnarray}
 I_{(a,b)}^{(\alpha)}(n,x)=\frac {1}{n^{\frac b a}}
 \left((\psi_1(x)
 +\frac {1}{12})x^{\alpha+\frac 1 a} -\frac {1}{12}n^{(a+b)(\alpha+\frac 1 a)}
 \right)
\end{eqnarray}
\begin{eqnarray*}
&&+\frac {1}{n^{\frac b a}}\times O\left((x+n^{a+b})^{(\alpha
+\frac
1 a-1)}\right)\\
 &&+\sum_{j=n}^{\left[\left(\frac
{x}{n^b}\right)^{\frac {1}{a}}\right]}
 \left(
\frac {1}{\alpha +1}\left(\zeta (-1-\alpha)-\frac
{\alpha}{2(2+\alpha)}\right)
 +\frac {1}{12}(j^an^b)^\alpha
 +O((j^an^b)^{(\alpha -1)})\right)\\
 &&+(n-\frac 1 2)\left(
\frac {1}{\alpha +1}\left(\zeta (-1-\alpha)-\frac
{\alpha}{2(2+\alpha)}\right)
 +\frac {1}{12}(n^{a+b})^\alpha
 +O((n^{a+b})^{(\alpha -1)})\right)\\
 &&-\left( \left(\frac {x}{n^b}\right)^{\frac
{1}{a}}-\psi\left( \left(\frac {x}{n^b}\right)^{\frac {1}{a}}
\right)\right)\left( \frac {1}{\alpha +1}\left(\zeta
(-1-\alpha)-\frac {\alpha}{2(2+\alpha)}\right)
 +(\psi_1(x)+\frac {1}{12})x^\alpha
 +O(x^{(\alpha -1)})\right)\\
 =&&-\frac {1}{24}n^{(a+b)\alpha}+
\psi\left( \left(\frac {x}{n^b}\right)^{\frac {1}{a}}
\right)\left(\psi_1(x)+\frac {1}{12}\right))x^\alpha
 +\sum_{j=n}^{\left[\left(\frac
{x}{n^b}\right)^{\frac {1}{a}}\right]}
 \left(
\frac {1}{12}(j^an^b)^\alpha
 +O((j^an^b)^{(\alpha -1)})\right)\\
  &&+O\left(\frac {x^{\alpha +\frac 1 a-1}}{n^{\frac b a}}+n^{(a+b)(\alpha-1)+1})\right).
\end{eqnarray*}
We write
\begin{eqnarray}
\sum_{j=n}^{\left[\left(\frac {x}{n^b}\right)^{\frac
{1}{a}}\right]}
 \left(
\frac {1}{12}(j^an^b)^\alpha
 +O((j^an^b)^{(\alpha -1)})\right)
 \end{eqnarray}
\begin{eqnarray*}
 =\frac {1}{12}n^{b\alpha}\sum_{j=n}^{\left[\left(\frac {x}{n^b}\right)^{\frac
{1}{a}}\right]}j^{a\alpha}+O\left(n^{b(\alpha
-1)}\sum_{j=n}^{\left[\left(\frac {x}{n^b}\right)^{\frac
{1}{a}}\right]}j^{a(\alpha-1)}\right).
\end{eqnarray*}
If $\alpha \not =-\frac 1 a$, by (2.19) we get
\begin{eqnarray}
\sum_{j=n }^{\left[\left(\frac{x}{n^b}\right)^{\frac
{1}{a}}\right]}j^{a\alpha} =\sum_{j\leq
\left(\frac{x}{n^b}\right)^{\frac{1}{a}}} j^{a\alpha} -\sum_{j\leq
n}j^{a\alpha}+n^{a\alpha}
\end{eqnarray}
\begin{eqnarray*}
 =&&\frac { \left(\frac
{x}{n^b}\right)^ {\frac {1}{a}(1+a\alpha) }-n^{1+a\alpha} }
{1+a\alpha} -\psi \left(\left( \frac {x}{n^b}\right)^{\frac
{1}{a}}\right) \left(\left(\frac {x}{n^b}\right)^{\frac
{1}{a}}\right)^{a\alpha}\\
&&+\psi(n)n^{a\alpha}+n^{a\alpha} +O \left(\left( \frac
{x}{n^b}\right)^{\frac{1}{a}(a\alpha-1)} +n^{a\alpha-1}
\right)\\
=&&\frac { \left(\frac {x}{n^b}\right)^ {\alpha+\frac {1}{a}
}-n^{1+a\alpha} } {1+a\alpha} -\psi \left(\left( \frac
{x}{n^b}\right)^{\frac {1}{a}}\right) \left(\frac
{x}{n^b}\right)^{\alpha}\\
&&+\frac 1 2n^{a\alpha} +O \left(\left( \frac
{x}{n^b}\right)^{\alpha-\frac{1}{a}} +n^{a\alpha-1} \right).
\end{eqnarray*}

If $\alpha  =-\frac 1 a$, by (2.20) we have
\begin{eqnarray}
\sum_{j=n }^{\left[\left(\frac{x}{n^b}\right)^{\frac
{1}{a}}\right]}j^{a\alpha} =\sum_{j\leq
\left(\frac{x}{n^b}\right)^{\frac{1}{a}}} j^{-1} -\sum_{j\leq
n}j^{-1}+n^{-1}
\end{eqnarray}
\begin{eqnarray*}
=\frac 1 a\left(\log x-(a +b)\log n\right) -\psi \left(\left(
\frac {x}{n^b}\right)^{\frac {1}{a}}\right) \left(\frac
{x}{n^b}\right)^{-\frac 1 a}+\frac {1}{2n} +O( n^{-2}).
\end{eqnarray*}

From (4.9) and (4.10), we also have
\begin{eqnarray}
\sum_{j=n }^{\left[\left(\frac{x}{n^b}\right)^{\frac
{1}{a}}\right]}j^{a(\alpha-1)}\ll \left\{
\begin{array}{ll}
\left(\frac {x}{n^b}\right)^ {\alpha+\frac {1}{a}
-1}+n^{1+a(\alpha-1)},
&\mbox{if $\alpha \not=1-\frac {1}{a}$ } \\
\log \frac{x}{n^{a+b}}+\frac 1 n,&\mbox{if $\alpha =1-\frac
{1}{a}$}
\end{array}
\right.
\end{eqnarray}
\begin{eqnarray*}
 \ll \left(\frac {x}{n^b}\right)^ {\alpha+\frac {1}{a}
-1}\log x+n^{1+a(\alpha-1)}\log x.
\end{eqnarray*}

 Now (4.5) follows from  (4.7),(4.8),(4.9) and (4.11),  (4.6) follows from  (4.7),(4.8),(4.10) and (4.11). This completes the proof of Lemma 4.2.

\section{\bf The proofs of Theorem 1 and Theorem 2}

We first prove Theorem 1. We take $f(n)=d(a,b;n)$ ,
$g(x)=\zeta(\frac b a)x^{\frac 1 a}+\zeta(\frac a b)x^{\frac 1
   b}$ in Lemma 2.1.
By Lemma 2.1 with $k=2$, we
have
\begin{eqnarray}
\sum_{n\leq x}\Delta^2(a,b;n)=
\left(\frac{1}{2}-\psi (x)
\right)\Delta^2(a,b;x)+\int_1^x \Delta^2(a,b;t)dt
\end{eqnarray}
\begin{eqnarray*}
 &&+ 2\int_1^x
\left(\frac{1}{2}-\psi (t) \right ) \left(\frac 1 a\zeta(\frac b
a)t^{\frac 1 a-1}+\frac 1 b\zeta(\frac a b)t^{\frac 1
   b-1}\right)\Delta (a,b;t)dt\\
=&& \left(\frac{1}{2}-\psi (x)
\right)\Delta^2(a,b;x)+\int_1^x \Delta^2(a,b;t)dt +T_1-2T_2,
\end{eqnarray*}
where
\begin{eqnarray}
T_1:=\int_1^x  \left(\frac 1 a\zeta(\frac b a)t^{\frac 1
a-1}+\frac 1 b\zeta(\frac a b)t^{\frac 1 b-1}\right)\Delta
(a,b;t)dt,
\end{eqnarray}
and
\begin{eqnarray}
T_2:=\int_1^x \left(\frac 1 a\zeta(\frac b a)t^{\frac 1 a-1}+\frac
1 b\zeta(\frac a b)t^{\frac 1 b-1}\right)\psi (t)\Delta (a,b;t)dt.
\end{eqnarray}

We treat $T_1$ first and shall  show that
\begin{eqnarray}
T_1=\frac {1} {4}\zeta(\frac b a )x^{\frac 1 a}+\frac {1}
  {4}\zeta(\frac a b)x^{\frac 1 b}
\end{eqnarray}
\begin{eqnarray*}
+\left(\frac 1 a\zeta(\frac b a)x^{\frac 1 a-1} +\frac 1
b\zeta(\frac a b)x^{\frac 1 b-1}\right) G_{(a,b)}(x) +O\left(
x^{\frac 1 a-\frac {3}{2(a+b)}}\right),
\end{eqnarray*}
where the series $G_{(a,b)}(x)$ was defined in   Theorem 1.

Integrating by parts, we have
\begin{eqnarray}T_1=
\left(\frac 1 a\zeta(\frac b a)x^{\frac 1 a-1}+\frac 1
b\zeta(\frac a b)x^{\frac 1 b-1}\right)
\int_1^x\Delta (a,b;t)dt
\end{eqnarray}
\begin{eqnarray*}
-\int_1^x \left(\frac 1 a(\frac 1 a-1)\zeta(\frac b a)t^{\frac 1
a-2}+\frac 1 b(\frac 1 b-1)\zeta(\frac a b)t^{\frac 1
b-2}\right)\left(\int_1^t\Delta (a,b;u)du\right)dt.
\end{eqnarray*}

We consider two cases.

\textbf{Case( i)}. If $a=1$ and $b\geq 2$, by Lemma 2.2 with $q=1$, a simple splitting argument and
 the first derivative test(See (2.3) in Ivi\'c[13])(Similar to the estimate of $T_1^*$ below in this paper),  we easily get
\begin{eqnarray*}T_1=&&\left(\zeta( b )+\frac {1} {b}\zeta(\frac 1 b)x^{\frac 1 b-1}\right)
\left(\frac {1} {4}x+G_{(a,b)}(x)+O\left( x^{1 -\frac {3}{2(1+b)}}\right)\right)\\
&&-\frac {1}{4b}(\frac 1 b-1)\zeta(\frac a b)b x^{\frac 1
b}+O\left( x^{\frac 1 b-\frac {3}{2(1+b)}}\right)\\
=&&\frac {1} {4}\zeta( b )x+\frac {1} {4}\zeta(\frac 1 b)x^{\frac
1 b}+\left(\zeta( b )+\frac {1} {b}\zeta(\frac 1 b)x^{\frac 1 b-1}\right)G_{(a,b)}(x)+O\left( x^{1-\frac {3}{2(1+b)}}\right).
\end{eqnarray*}
This proves that (5.4) holds in this case.

\textbf{Case (ii)}. Suppose $a\geq 2$ and $b>a$.  Similar to proof of the case (i), we also
have
\begin{eqnarray*}T_1=&&
\left(\frac 1 a\zeta(\frac b a)x^{\frac 1 a-1}+\frac 1
b\zeta(\frac a b)x^{\frac 1 b-1}\right)
\left(\frac {1} {4}x+G_{(a,b)}(x)+O\left( x^{ 1  -\frac {3}{2(a+b)}}\right)\right)\\
&&-\frac {1}{4a}(\frac 1 a-1)\zeta(\frac b a)a x^{\frac 1 a}-\frac
{1}{4b}(\frac 1 b-1)\zeta(\frac a b)b x^{\frac 1
b}+O\left( x^{\frac 1 a-\frac {3}{2(a+b)}}\right)\\
=&&\frac {1} {4}\zeta(\frac b a )x^{\frac 1 a}+\frac {1}
{4}\zeta(\frac a b)x^{\frac 1 b}+
 \left(\frac 1 a\zeta(\frac b a)x^{\frac 1 a-1}+\frac 1
 b\zeta(\frac a b)x^{\frac 1 b-1}\right)G_{(a,b)}(x)
+O\left( x^{\frac 1 a-\frac
{3}{2(a+b)}}\right).
\end{eqnarray*}
This completes the proof of (5.4).

Next we estimate $T_2$. From Lemma 3.1, we write
\begin{eqnarray*}
T_2=&&
-\int_1^x \left(\frac 1 a\zeta(\frac b a)t^{\frac 1
a-1}+\frac 1 b\zeta(\frac a b)t^{\frac 1 b-1}\right)\psi (t)
\left(\sum_{n^{a+b}\leq t}
\psi\left(\left(\frac{t}{n^{b}}\right)^{\frac
{1}{a}}\right)+\psi\left(\left(\frac{t}{n^{a}}\right)^{\frac {1}{
b}}\right)\right)dt\\
&&+\int_1^x \left(\frac 1 a\zeta(\frac b a)t^{\frac 1 a-1}+\frac 1
b\zeta(\frac a b)t^{\frac 1 b-1}\right)\psi (t)R(a,b;t)dt
\end{eqnarray*}
\begin{eqnarray}
= T_{21}+T_{22}.
\end{eqnarray}

To treat $T_{22}$, we divide the   interval $[1,x]$ into two subsets
$I_1$ and $I_2$, where $I_1=[1,x]\setminus\bigcup_{j=1}^{[x^{\frac
{1}{a+b}}]+1}[j^{a+b}-\frac {1}{10},j^{a+b}+\frac {1}{10}]$, and
$I_2=[1,x]\bigcap\left(\bigcup_{j=1}^{[x^{\frac
{1}{a+b}}]+1}[j^{a+b}-\frac {1}{10},j^{a+b}+\frac {1}{10}]\right)$.
 For $I_2$, we have trivial estimate
\begin{eqnarray*}
&&\int_{I_2} \left(\frac 1 a\zeta(\frac b a)t^{\frac 1 a-1}+\frac
1 b\zeta(\frac a b)t^{\frac 1 b-1}\right)\psi (t)R(a,b;t)dt\\
&&\ll \int_{I_2}t^{\frac 1 a-1}dt \ll \sum_{j\leq x^{\frac {1}{
a+b}}}j^{\frac 1 a-1}\ll x^{\frac {1}{ a(a+b)}}.
\end{eqnarray*}

For $I_1$, by applying integration by parts and (3.5) in Lemma
3.1, we get
\begin{eqnarray*} &&\int_{I_1} \left(\frac 1 a\zeta(\frac b
a)t^{\frac 1 a-1}+\frac
1 b\zeta(\frac a b)t^{\frac 1 b-1}\right)\psi (t)R(a,b;t)dt\\
&&=\sum_{j\leq x^{\frac {1}{
a+b}}}\int_{[1,x]\bigcap[j^{a+b}+\frac {1}{10},(j+1)^{a+b}-\frac
{1}{10}]}\left(\frac 1 a\zeta(\frac b a)t^{\frac 1 a-1}+\frac 1
b\zeta(\frac a b)t^{\frac 1 b-1}\right)\psi (t)R(a,b;t)dt\\
&&\ll\sum_{j\leq x^{\frac {1}{ a+b}}}\left( ((j+1)^{a+b})^{(\frac
1 a-1)}+\int_{[1,x]\bigcap[j^{a+b}+\frac {1}{10},(j+1)^{a+b}-\frac
{1}{10}]}|\psi_1(t)|\left(t^{\frac 1 a-2}+t^{\frac 1
a-1}|R'(a,b;t)|\right)dt
\right)\\
&&\ll\sum_{j\leq x^{\frac {1}{ a+b}}}\left(((j+1)^{a+b})^{(\frac 1
a-1)}+\int_{j^{a+b}}^{(j+1)^{a+b}}t^{\frac 1 a-1}t^{\frac
{1}{a+b}-1}dt\right)\\
&&\ll \sum_{j\leq x^{\frac {1}{ a+b}}}\left( ((j+1)^{a+b})^{(\frac
1 a-1)}+\left|(j+1)^{(a+b)(\frac 1 a+\frac
{1}{a+b}-1)}-j^{(a+b)(\frac 1 a+\frac
{1}{a+b}-1)}\right|\right) \\
&&\ll \sum_{j\leq x^{\frac {1}{ a+b}}}\left( ((j+1)^{a+b})^{(\frac
1 a-1)}+(j+1)^{(a+b)(\frac 1 a+\frac {1}{a+b}-1)-1}+j^{(a+b)(\frac
1 a+\frac
{1}{a+b}-1)-1}\right) \\
&&\ll \sum_{j\leq x^{\frac {1}{ a+b}}}\left( (j^{(a+b)})^{(\frac 1
a-1)}+(j+1)^{(a+b)(\frac 1 a-1)}\right) \\
&&\ll \left\{
\begin{array}{ll}
 x^{\frac{1}{a+b}}, &\mbox{if  $a=1$}\\
1,&\mbox{if $a\geq 2$}.
\end{array}
\right.
\end{eqnarray*}

Combining the above two estimates, we get
\begin{eqnarray}
T_{22}\ll x^{\frac{1}{a(a+b)}}.
\end{eqnarray}

To estimate $T_{21}$, we write

\begin{eqnarray*}
T_{21}&&=-\int_1^x \left(\frac 1 a\zeta(\frac b a)t^{\frac 1
a-1}+\frac 1 b\zeta(\frac a b)t^{\frac 1 b-1}\right)\psi (t)
\left(\sum_{n^{a+b}\leq t}
\psi\left(\left(\frac{t}{n^{b}}\right)^{\frac
{1}{a}}\right)+\psi\left(\left(\frac{t}{n^{a}}\right)^{\frac {1}{
b}}\right)\right)dt\\
&&=-\sum_{n^{a+b}\leq x}\int_{n^{a+b}}^x \left(\frac 1
a\zeta(\frac b a)t^{\frac 1 a-1}+\frac 1 b\zeta(\frac a b)t^{\frac
1 b-1}\right)\psi (t) \left(
\psi\left(\left(\frac{t}{n^{b}}\right)^{\frac
{1}{a}}\right)+\psi\left(\left(\frac{t}{n^{a}}\right)^{\frac {1}{
b}}\right)\right)dt\\
&&=-\sum_{n^{a+b}\leq x}\frac 1 a\zeta(\frac b
a)\left(I_{(a,b)}^{(\frac 1 a -1)}(n,x)+I_{(b,a)}^{(\frac 1 a
-1)}(n,x)\right)+\frac 1 b\zeta(\frac a b)\left(I_{(a,b)}^{(\frac 1
b-1)}(n,x)+I_{(b,a)}^{(\frac 1 b -1)}(n,x)\right)
\end{eqnarray*}
\begin{eqnarray}
:=-\frac 1 a\zeta(\frac b a)(T_{211}+T_{212})-\frac 1
b\zeta(\frac a b)(T_{213}+T_{214}),
\end{eqnarray}
and will show
\begin{eqnarray}
T_{21}=\left\{
\begin{array}{ll}
-\frac {1}{12}\zeta(b)\left(\zeta (b)x+2\zeta (\frac 1 b)x^{\frac 1 b}\right)
+O\left( x^{\frac {1}{1+b}}\right),&\mbox{if  $a=1$, and $b\geq 2$}\\
O(1+\log x),&\mbox{if $a\geq 2$ and $b>a$}.
\end{array}
\right.
\end{eqnarray}

From (2.19) and (2.20), we easily obtain
\begin{eqnarray}
\sum_{n\leq x^{\frac{1}{a+b}}} n^{-\frac
 b a}=\left\{
\begin{array}{ll}
\zeta(\frac b a) +O\left( x^{\frac {2}{a+b}-\frac 1 a}\right),&\mbox{if  $b>a,$}\\
\frac {a}{a-b}x^{\frac {2}{a+b}-\frac 1 a}+O(1),&\mbox{if $b<a$},
\end{array}
\right.
\end{eqnarray}
and
\begin{eqnarray}
\sum_{n\leq x^{\frac{1}{a+b}}} n^{(a+b)\alpha +1}=\left\{
\begin{array}{ll}
\frac {x^{\alpha+\frac{2}{a+b}}}{(a+b)\alpha +2}+O(1+x^{\alpha
+\frac {1}{a+b}}), &\mbox{if  $\alpha >-\frac {2}{a+b},$}\\
\frac {1}{a+b}\log x +\gamma +O(x^{-\frac {1}{a+b}}),&\mbox{if
$\alpha =-\frac {2}{a+b},$}\\
\zeta \left(-(a+b)\alpha -1\right)+O(x^{\alpha +\frac {2}{a+b}}),
&\mbox{if $\alpha <-\frac {2}{a+b}$}.
\end{array}
\right.
\end{eqnarray}

We will also use the following estimate
\begin{eqnarray}
\sum_{n\leq x}\frac{\log n}{n^{\frac 1 2}}=2(\log x-2)x^{\frac 1
2}+O(1),
\end{eqnarray}
which is an easy consequence of the  Euler-Maclaurin formula.

 Now we return to prove (5.9), and consider three cases.

\textbf{Case (1)}. $a=1$, $b=2$.

By (4.5) in Lemma 4.2, (5.10) and (5.11), we get
\begin{eqnarray*}
T_{211}=&&\sum_{n\leq x^{\frac 1 3}}I_{(1,2)}^{(0)}(n,x)\\
=&&\sum_{n\leq x^{\frac 1 3}}\frac {1} {12}\left(
\frac{x}{n^2}-n\right) + \psi\left( \frac {x}{n^2}\right)
\psi_1(x)\\
 &&+O\left(\sum_{n\leq x^{\frac 1 3}} \frac{1}{n^2}\log x+x^{-1}n^2+\frac 1 n
+\frac {1}{n^2}\log x\right)\\
=&&\frac {x}{12}\sum_{n\leq x^{\frac 1 3}}\frac{1}{n^2} -\frac
{1}{12}\left(\frac 1 2x^{\frac 2 3}-\psi (x^{\frac 1 3})x^{\frac 1
3}\right)+
 \psi_1(x)\sum_{n\leq x^{\frac 1 3}}
\psi\left( \frac {x}{n^2}\right) +O(\log x).
\end{eqnarray*}
(Here we use $\sum_{n\leq x^{\frac 1 3}}n=\frac 1 2x^{\frac 2 3}-\psi (x^{\frac 1 3})x^{\frac 1 3}+O(1)$).
  By (2.19) again, we obtain

\begin{eqnarray}
T_{211}=\frac {1}{12}\zeta(2)x-\frac {1}{8}x^{\frac 2 3}+O(x^{\frac 13}).
\end{eqnarray}

By (4.5) of Lemma 4.2, (5.10), (5.11) and (2.19), we have
\begin{eqnarray*}
T_{212}=&&\sum_{n\leq x^{\frac 1 3}}I_{(2,1)}^{(0)}(n,x)\\
=&&\sum_{n\leq x^{\frac 1 3}}\frac {1} {12}\left(
\frac{x^{\frac 1 2}}{n^{\frac 1 2}}-n\right) + \psi\left(\left( \frac {x}{n}\right)^{\frac 1 2}\right)
\psi_1(x)\\
 &&+O\left(\sum_{n\leq x^{\frac 1 3}} \frac{x^{-\frac 1 2}}{n^{\frac 1 2}}(1+\log x)+x^{-\frac 1 2}n^{\frac 1 2}+\frac 1 n
+\frac {1}{n^2}\log x\right)\\
=&&\frac {1}{12}\sum_{n\leq x^{\frac 1 3}}\frac{x^{\frac 1
2}}{n^{\frac 1 2}}-\frac {1}{12}\left(\frac 1 2x^{\frac 2 3}-\psi
(x^{\frac 1 3})x^{\frac 1 3}\right)
 +\psi_1(x)\sum_{n\leq
x^{\frac 1 3}} \psi\left(\left( \frac {x}{n}\right)^{\frac 1
2}\right)
+O(\log x)\\
=&&\frac 1 8x^{\frac 2 3}+\frac {1}{12}\zeta(\frac 1 2)x^{\frac 1 2}+O(x^{\frac 13}).
\end{eqnarray*}
In the same way, by (4.5) in Lemma 4.2, (5.10), (5.11) and (2.19),
we have
\begin{eqnarray*}
T_{213}=&&\sum_{n\leq x^{\frac 1 3}}I_{(1,2)}^{(-\frac 1 2)}(n,x)\\
=&&\frac {1} {6}\sum_{n\leq x^{\frac 1 3}}
\left(\frac{x^{\frac 1 2}}{n^2}-n^{-\frac 1 2} \right)+ \psi_1(x)x ^{ - \frac1 2}\sum_{n\leq x^{\frac 1 3}}\psi\left( \frac {x}{n^2}\right)
+O(\log x)\\
=&&\frac 1 6\zeta(2)x^{\frac 1 2}+O(x^{\frac 1 3}).
\end{eqnarray*}
 By (4.6) in Lemma 4.2, (5.10), (5.11) and (5.12), we have
\begin{eqnarray*}
T_{214}=&&\sum_{n\leq x^{\frac 1 3}}I_{(2,1)}^{(-\frac 1 2)}(n,x)\\
=&&\frac {1} {24}\sum_{n\leq x^{\frac 1 3}} \frac{\log x-3\log
n}{n^{\frac 1 2}} + \psi_1(x)x ^{ -\frac 1 2}\sum_{n\leq x^{\frac
1 3}}\psi\left(\left( \frac {x}{n}\right)^{\frac 1 2}\right)
+O(\log x)\\
=&&\frac 1 2 x^{\frac 1 6}++O(\log x).
\end{eqnarray*}
Combining the above four estimates , we see that (5.9) holds in
this case.

\textbf{Case (2)}. $a=1$, $b\geq 3$.

 By (4.5) in Lemma 4.2, (5.10), (5.11) and (2.19), we can get the following four
 estimates:
\begin{eqnarray*}
T_{211}=&&\sum_{n\leq x^{\frac {1}{1+b}}}I_{(1,b)}^{(0)}(n,x)\\
=&&\frac {1} {12}\sum_{n\leq x^{\frac {1}{1+b}}}\left(
\frac{x}{n^ b}-n\right) + \psi_1(x)\sum_{n\leq x^{\frac {1}{1+b}}}\psi \left(\frac
{x}{n^b}\right)+O(\log x)\\
=&&\frac {1}{12}\zeta(b)x-\frac {1}{12}(\frac {1}{b-1}+\frac 1 2)x^{\frac {2}{1+b}}+O(x^{\frac {1}{1+b}}),
\end{eqnarray*}

\begin{eqnarray*}
T_{212}=&&\sum_{n\leq x^{\frac {1}{1+b}}}I_{(b,1)}^{(0)}(n,x)\\
=&&\frac {1} {12}\sum_{n\leq x^{\frac {1}{1+b}}}\left(
\frac{x^{\frac 1 b}}{n^{\frac 1 b}}-n\right) + \psi_1(x)\sum_{n\leq x^{\frac {1}{1+b}}}\psi\left( \left(\frac
{x}{n}\right)^{\frac {1}{b}} \right)+O(\log x)\\
=&&\frac {1}{12}\left(\frac {1}{b-1}+\frac 1 2\right)x^{\frac {2}{b+1}}
+\frac {1}{12}\zeta (\frac 1 b)x^{\frac 1 b}+O(x^{\frac {1}{b+1}}),
\end{eqnarray*}

\begin{eqnarray*}
T_{213}=&&\sum_{n\leq x^{\frac {1}{1+b}}}I_{(1,b)}^{(-1+\frac 1
b)}(n,x)\\
=&&\frac {b} {12}\sum_{n\leq x^{\frac {1}{1+b}}}\left(
\frac{x^{\frac 1 b}}{n^ b
}-n^{-b+1+\frac 1 b}\right) + \psi_1(x)x^{-1+\frac 1 b}\sum_{n\leq x^{\frac {1}{1+b}}}\psi\left(\frac
{x}{n^b}\right)+O(\log x)\\
=&&\frac {b}{12}\zeta(b)x^{\frac 1 b}+O(\log x),
\end{eqnarray*}
and
\begin{eqnarray*} T_{214}=&&\sum_{n\leq x^{\frac
{1}{1+b}}}I_{(b,1)}^{(-1+\frac 1
b)}(n,x)\\
=&&\frac {1} {12(2-b)}\sum_{n\leq x^{\frac {1}{1+b}}}\left(
\frac{x^{-1 +\frac 2 b}}{n^{\frac 1 b}}-n^{-b+1+\frac 1 b}\right)
+ \psi_1(x)x^{-1+\frac 1 b}\sum_{n\leq x^{\frac {1}{1+b}}}\psi\left( \left(\frac
{x}{n}\right)^{\frac {1}{b}} \right)+O(\log x)\\
=&&O(\log x).
\end{eqnarray*}
Hence, (5.9) also holds in this case.

\textbf{Case (3)}. $a\geq 2$, $b\geq 3$.

 Similar to the proof of Case (2), we can prove

\begin{eqnarray}
T_{211},T_{212},T_{213},T_{214}\ll  \log x.
\end{eqnarray}
We omit the details. This completes the the proof of (5.9).

Note that if $a=1, b\geq 2$, then $\frac {1}{1+b}<1-\frac
{3}{2(1+b)}$; if $a\geq 2, b>a$, then $\frac {1}{a(a+b)}<\frac 1
a-\frac {3}{2(a+b)}$. Now, collecting (5.1), (5.5)-(5.9) completes
the proof of Theorem 1.

\bigskip

Finally, we shall give a short  proof of Theorem 2, since the
details are   similar to  and simpler than that of Theorem 1.

We take $f(n)=d(1,1;n)$ , $g(x)=(\log x +2\gamma -1)x$. By Lemma 2.1 with $k=2$, we have
\begin{eqnarray}
\sum_{n\leq x}\Delta^2(1,1;n)=
\left(\frac{1}{2}-\psi (x)
\right)\Delta^2(1,1;x)+\int_1^x \Delta^2(1,1;t)dt
+T_1^*+T_2^*,
\end{eqnarray}
where
\begin{eqnarray}
T_1^*:=\int_1^x  (\log t +2\gamma)\Delta
(1,1;t)dt,
\end{eqnarray}
and
\begin{eqnarray}
T_2^*:=-2\int_1^x (\log t +2\gamma )\psi (t)\Delta (1,1;t)dt.
\end{eqnarray}

Similar to the estimate of $T_2$, we may get that( Also see
Furuya[5],page 17-18)
\begin{eqnarray}
T_2^*:=\frac 1 6x\log^2x+\frac 1 3(2\gamma -1)x\log x+\frac 1 3(2\gamma^2-2\gamma +1)x+O(x^{\frac 1 2}\log x).
\end{eqnarray}
Now we estimate $T_1^*$ as $T_1$ . By Lemma 2.2 with $q=1$, and
then integrating by parts, we find that
\begin{eqnarray*}
T_1^*=&&\int_1^x  (\log t +2\gamma)\Delta(1,1;t)dt \\
 =&&(\log x+2\gamma)\int_1^x \Delta(1,1;t)dt-\int_1^x  t^{-1}\left(\int_1^t\Delta(1,1;u)du\right)dt
 \\
 =&&(\log x+2\gamma)\left(\frac 1 4 x+G_{(1,1)}(x)+O(x^{\frac 1 4})\right)-
 \int_1^x  t^{-1}\left(\frac t 4+G_{(1,1)}(t)+O(t^{\frac 1
 4})\right)dt\\
 =&& \frac 1 4 \left(\log x+2\gamma-1\right)x+(\log x+2\gamma)G_{(1,1)}(x)-
 \int_1^x  t^{-1}G_{(1,1)}(t)dt+O(x^{\frac 1 4}\log x),
\end{eqnarray*}
where
\begin{eqnarray*}
 G_{(1,1)}(x):=\frac {1}{2{\sqrt 2} \pi^2}x^{\frac 3 4}
 \sum_{n=1}^{\infty}\frac {d(n)}{n^{\frac 5 4}}
 \sin \left( 4\pi \sqrt {nx}-\frac {\pi}{4}\right).
 \end{eqnarray*}
Since the  series $G_{(1,1)}(x)$ is absolute convergent,  one may
integrate term by term, and obtain

\begin{eqnarray*}
 \int_1^x  t^{-1}G_{(1,1)}(t)dt=
 \frac {1}{2{\sqrt 2} \pi^2}
 \sum_{n=1}^{\infty}\frac {d(n)}{n^{\frac 5 4}}
 \int_1^xt^{-\frac 1 4}\sin \left( 4\pi \sqrt {nt}-\frac {\pi}{4}\right)dt.
\end{eqnarray*}

 By  a simple splitting argument and
 the first derivative test(see (2.3), Ivi\'c[13]),  we easily get

 \begin{eqnarray}
 \int_1^x  t^{-1}G_{(1,1)}(t)dt\ll
 \sum_{n=1}^{\infty}\frac {d(n)}{n^{\frac 5 4}}\frac {x^{\frac 1
 4}}{\sqrt n}\ll x^{\frac 1 4}.
 \end{eqnarray}

 Hence,
 \begin{eqnarray}
 T_1^*=\frac 1 4 \left(\log x+2\gamma-1\right)x+(\log
 x+2\gamma)G_{(1,1)}(x)+O(x^{\frac 1 4}).
 \end{eqnarray}

Combining (5.15), (5.18) and (5.20) completes the proof of Theorem
2.

Xiaodong Cao\\ Dept. of Mathematics and Physics,\\ Beijing Institute
of
Petro-Chemical Technology,\\
Beijing, 102617, P.R. China\\
Email: caoxiaodong@bipt.edu.cn

\bigskip

Wenguang Zhai,\\
School of Mathematical Sciences,\\
Shandong Normal University,\\Jinan,  Shandong, 250014,\\
P.R.China\\
 E-mail:  zhaiwg@hotmail.com

\end{document}